\documentclass[11pt]{article}

\usepackage{mathrsfs}
\usepackage{amssymb}

\newtheorem{lemma}{Lemma}[section]
\newtheorem{theorem}{Theorem}[section]
\newtheorem{corollary}{Corollary}[section]

\def\scr{\mathscr}
\def\mbb{\mathbb}

\def\blemma{\begin{lemma}\sl{}\def\elemma{\end{lemma}}}
\def\btheorem{\begin{theorem}\sl{}\def\etheorem{\end{theorem}}}

\def\beqlb{\begin{eqnarray}}\def\eeqlb{\end{eqnarray}}
\def\beqnn{\begin{eqnarray*}}\def\eeqnn{\end{eqnarray*}}

\def\qed{\quad$\Box$\medskip}

\def\<{\langle}\def\>{\rangle}

\def\E{ \mbox{\boldmath $E$} }
\def\P{ \mbox{\boldmath $P$} }

\def\itOmega{{\it\Omega}}

\def\IR{\mbb{R}}

\begin{document}

\

\noindent{Published in: {\it Potential Analysis} {\bf 21} (2004),
1: 75--97}

\bigskip\bigskip

\noindent{\Large\bf Generalized Mehler Semigroups and Catalytic}

\medskip

\noindent{\Large\bf Branching Processes with
Immigration\footnote{Supported by a Max Planck Award (D.A.D.), the
NSERC (D.A.D. and B.S.), the NNSF (Z.L.) and a MITACS-PINTS
Postdoctoral Fellowship (W.S.).}}

\bigskip
\noindent{DONALD A. DAWSON}

\smallskip
\noindent{School of Mathematics and Statistics, Carleton
University,}

\noindent{Ottawa, Ontario, Canada K1S 5B6}

\noindent{(e-mail: {\tt ddawson@math.carleton.ca})}

\bigskip

\noindent{ZENGHU LI}

\smallskip

\noindent{Department of Mathematics, Beijing Normal University,}

\noindent{Beijing 100875, P.R. China}

\noindent{(e-mail: {\tt lizh@email.bnu.edu.cn})}

\bigskip
\noindent{BYRON SCHMULAND and WEI SUN}

\smallskip

\noindent{Department of Mathematical and Statistical Sciences, University of
Alberta,}

\noindent{Edmonton, Alberta, Canada T6G 2G1}

\noindent{(e-mail: {\tt schmu@stat.ualberta.ca} and {\tt
wsun@math.ualberta.ca})}

\bigskip\bigskip

\noindent{\bf Abstract.} Skew convolution semigroups play an
important role in the study of generalized Mehler semigroups and
Ornstein-Uhlenbeck processes. We give a characterization for a
general skew convolution semigroup on real separable Hilbert space
whose characteristic functional is not necessarily differentiable
at the initial time. A connection between this subject and
catalytic branching superprocesses is established through
fluctuation limits, providing a rich class of non-differentiable
skew convolution semigroups. Path regularity of the corresponding
generalized Ornstein-Uhlenbeck processes in different topologies
is also discussed.

\bigskip\bigskip


\bigskip


\bigskip\bigskip

\section{Introduction}
\setcounter{equation}{0}

Suppose that $H$ is a real separable Hilbert space. Given a Borel
probability measure $\nu$ on $H$, let $\hat\nu$ denote its
characteristic functional. It is known that if $\nu$ is infinitely
divisible, then $\hat \nu(a) \neq 0$ for all $a\in H$ and there is
a unique continuous function $\log\hat \nu$ on $H$ such that
$\log\hat \nu(0) =0$ and $\hat \nu(a) = \exp\{\log\hat \nu(a)\}$;
see e.g.\ Linde \cite[p.20 and p.58]{L86}. Let $(T_t)_{t\ge0}$ be
a strongly continuous semigroup of linear operators on $H$ with
dual $(T_t^*)_{t\ge0}$ and $(\mu_t)_{t\ge0}$ a family of
probability measures on $H$. The family $(\mu_t)_{t\ge0}$ is
called a {\it skew convolution semigroup} (SC-semigroup)
associated with $(T_t)_{t\ge0}$ if the following equation is
satisfied:
 \beqlb\label{1.1}
\mu_{r+t} = (T_t\mu_r) * \mu_t, \qquad r,t\ge0,
 \eeqlb
where ``$*$'' denotes the convolution operation. It
is easy to check that (\ref{1.1}) holds if and only if we can
define a Markov transition semigroup $(Q^\mu_t)_{t\ge0}$ on $H$ by
 \beqlb\label{1.2}
Q^\mu_tf(x) := \int_H f(T_tx+y)\mu_t(dy), \quad x\in H,f\in B(H),
 \eeqlb
where $B(H)$ denotes the totality of bounded Borel measurable
functions on $H$. In this case, $(Q^\mu_t)_{t\ge0}$ is called a
{\it generalized Mehler semigroup}, which corresponds to a
generalized Ornstein-Uhlenbeck process (OU-process) with state
space $H$. This formulation of OU-processes was given by Bogachev
{\it et al} \cite{BRS96} as a generalization of the classical
Mehler formula; see e.g.\ Malliavin \cite[p.17 and p.25]{M97}. One
motivation to study such OU-processes is that they constitute a
large class of explicit examples of processes on
infinite-dimensional spaces with rich mathematical structures.
They arise in the study of Langevin type equations with
generalized drift involving the generator of $(T_t)_{t\ge0}$. We
refer the reader to Bogachev {\it et al} \cite{BRS96}, Fuhrman and
R\"ockner \cite{FR00}, and van Neerven \cite{N00} for discussions
from a theoretical viewpoint. See also Bogachev and R\"ockner
\cite{BR95}, Fuhrman \cite{F98}, and van Neerven \cite{N98} for
some earlier related work. In the setting of cylindrical
probability measures, Bogachev {\it et al} \cite[Lemma~2.6]{BRS96}
proved that, if the function $t \mapsto \hat\mu_t(a)$ is
absolutely continuous on $[0,\infty)$ and differentiable at $t=0$
for all $a\in H$, then (\ref{1.1}) is equivalent to
 \beqlb\label{1.3}
\hat \mu_t(a) = \exp\bigg\{-\int_0^t\lambda(T^*_sa)ds\bigg\},
\quad t\ge0, a\in H,
 \eeqlb
where $\lambda(a) = -(d/dt) \hat \mu_t(a)|_{t=0}$ is a
negative-definite functional on $H$. A necessary and sufficient
condition for a Gaussian SC-semigroup to be differentiable was
given in van Neerven \cite{N00}. These results give
characterizations for interesting special classes of SC-semigroups
defined by (\ref{1.1}) and have stimulated the present work.

Skew convolution semigroups have also played an important role in
the study of immigration structures associated with branching
processes. Let $E$ be a Lusin topological space, i.e., a
homeomorph of a Borel subset of a compact metric space, with Borel
$\sigma$-algebra ${\scr B}(E)$. We denote by $B(E)^+$ the set of
bounded non-negative Borel functions on $E$. Let $M(E)$ be the
totality of finite measures on $(E,{\scr B}(E))$ endowed with the
topology of weak convergence and $(Q_t)_{t\ge0}$ the transition
semigroup of a {\it measure-valued branching process}
(superprocess) $X$ with state space $M(E)$. A family
$(N_t)_{t\ge0}$ of probability measures on $M(E)$ is called a {\it
SC-semigroup} associated with $(Q_t)_{t\ge0}$ if it satisfies
 \beqlb\label{1.4}
N_{r+t} = (N_rQ_t)*N_t,
\qquad r,t\ge0.
 \eeqlb
We use the same terminology for solutions of (\ref{1.1}) and those
of (\ref{1.4}) since (\ref{1.1}) is actually a special form of
(\ref{1.4}) when they are put in a slightly more general setting,
say, when $H$ and $M(E)$ are replaced by a topological semigroup.
This similarity between the two equations was first noticed by
L.G. Gorostiza (1999, personal communication); see also Bojdecki
and Gorostiza \cite{BG02} and Schmuland and Sun \cite{SS01}. It is
not hard to show that (\ref{1.4}) holds if and only if
 \beqlb\label{1.5}
Q^N_t(\nu,\cdot):= Q_t(\nu,\cdot)*N_t,
\qquad t\ge0, \nu\in M(E)
 \eeqlb
defines a Markov semigroup $(Q^N_t)_{t\ge0}$ on $M(E)$. A Markov
process $Y$ in $M(E)$ is called an {\it immigration process}
associated with $X$ if it has transition semigroup $(Q^N_t)
_{t\ge0}$. The intuitive meaning of the immigration process is
clear from (\ref{1.5}), that is, $Q_t(\nu,\cdot)$ is the
distribution of descendants of the people distributed as $\nu\in
M(E)$ at time zero and $N_t$ is the distribution of descendants of
the people immigrating to $E$ during the time interval $(0,t]$. By
Li \cite[Theorem~2]{L95} or \cite[Theorem~3.2]{L02}, the family
$(N_t)_{t\ge0}$ satisfies (\ref{1.4}) if and only if there is an
infinitely divisible probability entrance law $(K_s)_{s>0}$ for
$(Q_t)_{t\ge0}$ such that
 \beqlb\label{1.6}
\log \int_{M(E)} e^{-\nu(f)} N_t(d\nu)
=
\int_0^t\bigg[\log\int_{M(E)} e^{-\nu(f)} K_s(d\nu)\bigg] d s ,
\quad t\ge0, f\in B(E)^+,
 \eeqlb
where $\nu(f) = \int_E f d\nu$; see also Li \cite{L98, L02} for
some generalizations of this result. Then there is a 1-1
correspondence between SC-semigroups and a set of infinitely
divisible probability entrance laws. Some representations of the
infinitely divisible probability entrance laws and path regularity
of the corresponding immigration processes were studied in Li
\cite{L96}. The connection between immigration processes and
generalized OU-processes was studied in Gorostiza and Li \cite{GL98,
GL00} and Li \cite{L99}. In view of (\ref{1.6}), the function
 \beqlb\label{1.7}
t\mapsto
\log\int_{M(E)} e^{-\nu(f)} N_t(d\nu)
 \eeqlb
is always absolutely continuous on $[0,\infty)$, and it is
differentiable at $t=0$ for all continuous $f\in B(E)^+$ if and
nearly only if $(K_s)_{s>0}$ is closable by an infinitely
divisible probability measure $K_0$ on $M(E)$. By the similarity
of (\ref{1.1}) and (\ref{1.4}), one might expect similar results
for the solutions of (\ref{1.1}). However, the Hilbert space
situation is much more complicated as Schmuland and Sun \cite{SS01}
showed that the linear part of $t\mapsto \log\hat \mu_t(a)$ can be
discontinuous. Therefore, we can only discuss characterizations
for the solutions of (\ref{1.1}) under reasonable regularity
conditions on the linear part of $t\mapsto \log\hat \mu_t(a)$.

This work is also related to the catalytic branching superprocess
introduced by Dawson and Fleischmann \cite{DF91, DF92}. Let us
consider the special case where the underlying motion is an
absorbing barrier Brownian motion (ABM) in a domain $D$. Let
$(P_t)_{t\ge0}$ denote the transition semigroup of the ABM. Let
$\eta \in M(D)$ and let $\phi (\cdot, \cdot)$ be a function on
$D\times [0,\infty)$ of a certain form to be specified. A {\it
catalytic branching superprocess} in $M(D)$ has transition
semigroup $(Q_t)_{t\ge0}$ determined by
 \beqlb\label{1.8}
\int_{M(D)} e^{-\nu(f)} Q_t(\mu,d\nu)
=
\exp\left\{-\mu(V_tf)\right\},
\qquad f\in B(D)^+,
 \eeqlb
where $(V_t)_{t\ge0}$ is a semigroup of non-linear operators on
$B(D)^+$ defined by
 \beqlb\label{1.9}
V_tf(x) = P_t f(x)
- \int_0^td s\int_D \phi(y,V_sf(y)) p_{t-s}(x,y)\eta(d y),
\quad t\ge0, x\in D,
 \eeqlb
with $p_t(x,y)$ being the density of $P_t(x,dy)$. This process
describes the catalytic reaction of a large number of
infinitesimal particles moving in $D$ according to the transition
law of the ABM and splitting according to the branching mechanism
given by $\phi(\cdot,\cdot)$. The measure $\eta(dx)$ represents
the distribution of a catalyst in $D$ which causes the splitting.
More detailed descriptions of the model will been given in Section 3.

In this paper, we give a representation for the general
SC-semigroup $(\mu_t)_{t\ge0}$ defined by (\ref{1.1}) whose
characteristic functional is not necessarily differentiable at
$t=0$. This result extends the interesting characterizations given
in \cite{BRS96} and \cite{N00}. The general representation is of
interest since it includes some SC-semigroups arising in
applications which are not included in (\ref{1.3}). We provide a
rich class of SC-semigroups of this type in the case where $H=
L^2(0,\infty)$ and $(T_t)_{t\ge0}$ is the transition semigroup of
the ABM. Indeed, the corresponding generalized OU-processes arise
naturally as fluctuation limits of catalytic branching
superprocesses with immigration. An important feature of these
OU-processes is that they usually do not have right continuous
realizations, which is similar to the situation of immigration
processes studied in \cite{L96, L98, L02}. Nevertheless, we show
that some of these OU-processes are in fact quite regular if we
regard them as processes with values of signed-measures. The study
of generalized Mehler semigroups on Hilbert spaces and that of
catalytic branching processes have evolved independently of
each other with different motivations, techniques, and so on. The
fluctuation limits establish a connection between the two
subjects.

The remainder of this paper is organized as follows. In Section 2
we give the characterization for general SC-semigroups.
Fluctuation limits of immigration processes are studied in Section
3, which lead to generalized OU-processes with distribution
values. Under stronger assumptions, it is proved in Section 4 that
some of these OU-processes actually live in the Hilbert space
$L^2(D)$ of functions. Regularity properties of the processes in
the space of signed-measures are discussed in Section 5.

\section{Characterization of SC-semigroups}
\setcounter{equation}{0}

In this section, we give a general representation of the
SC-semigroups defined by (\ref{1.1}). It was proved in Schmuland
and Sun \cite{SS01} that, if $(\mu_t)_{t\ge0}$ is a solution of
(\ref{1.1}), then each $\mu_t$ is an infinitely divisible probability
measure. Let
 \beqnn
K(x,a)
:=
e^{i\< x,a\>} -1 - i\< x,a\>\chi_{[0,1]}(\|x\|),
\qquad x, a\in H.
 \eeqnn
By Linde \cite[p.75 and p.84]{L86}, the characteristic functional
of $\mu_t$ on $H$ is given by
 \beqlb\label{2.1}
\hat\mu_t(a) = \exp\bigg\{i\<b_t,a\>
- \frac{1}{2}\<R_ta,a\> + \int_{H^\circ} K(x,a) M_t(dx)\bigg\},
\qquad a\in H,
 \eeqlb
where $b_t\in H$, $R_t$ is a symmetric, positive-definite, nuclear operator on $H$,
and $M_t$ is a $\sigma$-finite measure (L\'evy measure) on $H^\circ := H
\setminus \{0\}$ satisfying
 \beqlb\label{2.2}
\int_{H^\circ} (1\land\|x\|^2) M_t(dx)
<
\infty.
 \eeqlb
Thus, $\mu_t$ is uniquely determined by the triple $(b_t,R_t,M_t)$
and is uniquely decomposed into the convolution of three
infinitely divisible probabilities $\mu_t = \mu^c_t * \mu^g_t
* \mu^j_t$ with
 \beqlb\label{2.3}
\hat\mu^c_t(a) = \exp\{i\<b_t,a\>\},
\qquad
\hat\mu^g_t(a) = \exp\bigg\{- \frac{1}{2}\<R_ta,a\>\bigg\},
 \eeqlb
and
 \beqlb\label{2.4}
\hat\mu^j_t(a) = \exp\bigg\{\int_{H^\circ} K(x,a) M_t(dx)\bigg\}
 \eeqlb
for  $a\in H$. We call $\mu^c_t$ the {\it constant} (or {\it
linear}) {\it part}, $\mu^g_t$ the {\it Gaussian part}, and
$\mu^j_t$ the {\it jump part} of $\mu_t$. By the uniqueness of the
decomposition (\ref{2.1}) it is not hard to show that (\ref{1.1})
holds if and only if we have
 \beqlb\label{2.5}
R_{r+t}
=
T_tR_rT_t^* + R_t,
\qquad
M_{r+t}
=
(T_tM_r)|_{H^\circ} + M_t,
 \eeqlb
and
 \beqlb\label{2.6}
b_{r+t}
=
b_t + T_tb_r
+ \int_{H^\circ} \left(\chi_{[0,1]}(\|T_tx\|)
- \chi_{[0,1]}(\|x\|)\right)T_tx\, M_r(dx)
 \eeqlb
for all $r,t\ge0$.

\btheorem\label{t2.1} If $(\mu_t)_{t\ge0}$ is an SC-semigroup
with decomposition (\ref{2.1}), then we can write
 \beqlb\label{2.7}
\<R_ta,a\>
=
\int_0^t \< U_sa,a\>ds,
\qquad t\ge0,a\in H,
 \eeqlb
where $(U_s)_{s>0}$ is a family of nuclear operators on $H$
satisfying $U_{s+t} = T_tU_sT_t^*$ for all $s,t>0$ and
 \beqnn
\int_0^t \mbox{\rm Tr}\, U_s ds <\infty,
\qquad t\ge0.
 \eeqnn
\etheorem

The basic idea of the proof of this theorem is similar to that of
\cite[Theorem~2]{L95}, but the argument in the present case is
more involved. We first prove two lemmas.

\blemma\label{l2.1} Under the conditions of Theorem~\ref{t2.1}, the
function $t\mapsto \<R_ta,b\>$ is absolutely continuous in $t\ge0$
for all $a,b\in H$. \elemma

\noindent{\it Proof.} If $(\mu_t)_{t\ge0}$ is an SC-semigroup, so is
$(\mu^g_t)_{t\ge0}$ by the first equation in (\ref{2.5}). Then we have
 \beqlb\label{2.8}
\int_H \|x\|^2 \mu^g_{r+t}(dx)
=
\int_H \|T_tx\|^2 \mu^g_r(dx) + \int_H \|x\|^2 \mu^g_t(dx),
\qquad r,t\ge0.
 \eeqlb
It follows that
 \beqlb\label{2.9}
g(t) := \int_H \|x\|^2 \mu^g_t(dx),
\qquad t\ge0
 \eeqlb
is a non-decreasing function. Since $(T_t)_{t\ge0}$ is strongly
continuous, there are constants $c\ge 1$ and $b\ge0$ such that
$\|T_t\| \le c e^{bt}$. We claim that, for $0 <r_1 <t_1 <\cdots
<r_n <t_n \le l$,
 \beqlb\label{2.10}
\sum_{j=1}^n[g(t_j) - g(r_j)]
\le
c^2 e^{2bl} g(\sigma_n),
 \eeqlb
where $\sigma_n = \sum_{j=1}^n (t_j-r_j)$. When $n=1$, this follows
from (\ref{2.8}). Now assume that (\ref{2.10}) holds for $n-1$. Applying
(\ref{2.8}) twice,
 \beqnn
\sum_{j=1}^n[g(t_j) - g(r_j)]
&\le&
[g(t_n) - g(r_n)]
+ c^2e^{2bl} \int_H \|x\|^2 \mu^g_{\sigma_{n-1}}(dx)   \\
&=&
\int_H \|T_{r_n}x\|^2 \mu^g_{t_n-r_n}(dx)
+ c^2e^{2bl} \int_H \|x\|^2 \mu^g_{\sigma_{n-1}}(dx)   \\
&\le&
c^2e^{2bl} \int_H \|T_{\sigma_{n-1}}x\|^2 \mu^g_{t_n-r_n}(dx)
+ c^2e^{2bl} \int_H \|x\|^2 \mu^g_{\sigma_{n-1}}(dx)  \\
&=&
c^2e^{2bl} \int_H \|x\|^2 \mu^g_{\sigma_n}(dx),
 \eeqnn
which gives (\ref{2.10}). Letting $r\to 0$ and $t\to 0$ in (\ref{2.8})
and using the fact that $g$ is a non-decreasing function
one sees that $g(t) \to0$ as $t\to 0$. By this and (\ref{2.10}), $g$
is absolutely continuous in $t\ge0$. From (\ref{2.5}) we see that
$\<R_ta,a\>$ is a non-decreasing function of $t\ge0$ for any $a\in
H$. For $t\ge r\ge 0$, (\ref{2.5}) yields
 \beqnn
\<R_ta,a\> - \<R_ra,a\>
&=&
\<R_{t-r}T_r^*a,T_r^*a\>
=
\int_H\<x,T_r^*a\>^2\mu^g_{t-r}(dx)   \\
&\le&
\|a\|^2\int_H \|T_rx\|^2 \mu^g_{t-r}(dx)
=
\|a\|^2[g(t) - g(r)].
 \eeqnn
Then $\<R_ta,a\>$ is absolutely continuous in $t\ge0$. Polarization shows
 that $\<R_ta,b\>$ is absolutely continuous in
$t\ge0$ for all $a,b\in H$. \qed

\blemma\label{l2.2} Under the condition of Theorem~\ref{t2.1},
there is a family of nuclear operators $(U_s)_{s>0}$ on $H$ such
that (\ref{2.7}) holds. \elemma

\noindent{\it Proof.} Let $\{e_n: n=1,2,\dots\}$ be an orthonormal
basis of $H$. By Lemma~\ref{l2.1}, there are locally integrable
functions $A_{m,n}$ on $[0,\infty)$ such that
 \beqlb\label{2.11}
\<R_te_m,e_n\>
=
\int_0^t A_{m,n}(s)ds,
\qquad  t\ge0,\ m,n\ge 1.
 \eeqlb
From the symmetry of $R_t$ we get
 \beqlb\label{2.12}
\int_0^t A_{m,n}(s)ds
=
\int_0^t A_{n,m}(s)ds,
 \eeqlb
while the positivity of $R_t$ gives
 \beqlb\label{2.13}
\<R_ta,a\>
=
\int_0^t \sum_{m,n=1}^\infty A_{m,n}(s)\<a,e_m\>\<a,e_n\>ds
\ge 0
 \eeqlb
for $a\in \mbox{span} \{e_1,e_2,\dots\}$. (The sum is actually
finite!) In addition, since $R_t$ is nuclear we have
 \beqlb\label{2.14}
\int_0^t \bigg(\sum_{n=1}^\infty A_{n,n}(s)\bigg) ds
=
\sum_{n=1}^\infty \<R_te_n,e_n\>
=
\mbox{Tr}(R_t)
<\infty.
 \eeqlb
Let $F$ be the Borel subset of $[0,\infty)$ consisting of all $s\ge 0$
such that $A_{m,n}(s) = A_{n,m}(s)$ for $m,n\ge1$ and
 \beqlb\label{2.15}
\sum_{m,n=1}^\infty A_{m,n}(s)\<a,e_m\>\<a,e_n\>
\ge 0
\quad
\mbox{and}
\quad
\sum_{n=1}^\infty A_{n,n}(s)
<\infty
 \eeqlb
for $a\in \mbox{span} \{e_1,e_2, \dots\}$ with rational
coefficients. As observed in the proof of Lemma~\ref{l2.1},
$\<R_ta,a\>$ is a non-decreasing function of $t\ge0$. By
(\ref{2.12}), (\ref{2.13}) and (\ref{2.14}), $F$ has full Lebesgue
measure. For any $s\in F$,
 \beqlb\label{2.16}
U_s a
=
\sum_{m,n=1}^\infty A_{m,n}(s)\<a,e_m\> e_n,
 \eeqlb
defines a positive-definite, symmetric linear operator on $\mbox{span}
\{e_1,e_2,\dots\}$. Taking $b= xe_m + ye_n$, with $x,y$ rational, we get
 \beqnn
\<U_sb,b\>
=
(x \,\ y)
\left(\begin{array}{cc}
A_{m,m}(s) &A_{m,n}(s)  \\
A_{n,m}(s) &A_{n,n}(s)
\end{array}\right)
\left(\begin{array}{c}
x  \\
y
\end{array}\right)
\ge 0,
 \eeqnn
so that the $2\times 2$ matrix above is non-negative definite.
Therefore, its determinant is non-negative, that is,
 \beqlb\label{2.17}
A_{m,n}(s)^2
\le
A_{m,m}(s)A_{n,n}(s).
 \eeqlb
Combined with the Cauchy-Schwarz inequality this gives,
 \beqnn
\|U_sa\|^2
&=&
\sum_{n=1}^\infty \bigg(\sum_{m=1}^\infty A_{m,n}(s)
\<a,e_m\>\bigg)^2   \\
&\le&
\sum_{n=1}^\infty A_{n,n}(s) \bigg(\sum_{m=1}^\infty
 A_{m,m}(s)^{1/2} |\<a,e_m\>|\bigg)^2   \\
&\le&
\bigg(\sum_{n=1}^\infty A_{n,n}(s)\bigg)^2 \|a\|^2
 \eeqnn
for $s\in F$ and  $a\in \mbox{span} \{e_1,e_2, \dots\}$.
 This shows that $U_s$ is a bounded operator and can
be extended to the entire space $H$. In fact, $U_s$ is a nuclear
operator since
 \beqnn
\mbox{Tr}(U_s)
=
\sum_{n=1}^\infty \<U_se_n,e_n\>
=
\sum_{n=1}^\infty A_{n,n}(s)
< \infty.
 \eeqnn
By (\ref{2.11}) and (\ref{2.16}), for $a\in \mbox{span} \{e_1,e_2,
\dots\}$ we have
 \beqlb\label{2.18}
\<R_ta,a\>
=
\sum_{m,n=1}^\infty \<a,e_m\>\<a,e_n\> \<R_te_m,e_n\>
=
\int_0^t\<U_sa,a\>ds,
\qquad t\ge0.
 \eeqlb
Since $s\mapsto \mbox{Tr}(U_s)$ is locally integrable, by dominated
convergence we see that (\ref{2.18}) holds for all $a\in H$.
For $s\not\in F$, we let $U_s$ be the zero operator. \qed

\noindent{\it Proof of Theorem 2.1.} Let $(U_s)_{s>0}$ be provided by Lemma 2.2.
Note that (\ref{2.7}) and the
first equation of (\ref{2.5}) imply
 \beqnn
\int_0^r \<U_{s+t}a,a\> ds
=
\int_0^r \<U_sT_t^*a,T_t^*a\> ds,
\qquad r,t\ge0, a\in H.
 \eeqnn
Since $H$ is separable, by Fubini's theorem, there are subsets
$G$ and $G_s$ of $[0,\infty)$ with full Lebesgue measure such that
 \beqnn
U_{s+t} = T_tU_sT_t^*,
\qquad s\in G, t\in G_s.
 \eeqnn
Choose a decreasing sequence $s_n\in G$ with $s_n\to 0$, and define
 \beqnn
\widetilde U_t := T_{t-s_n}U_{s_n}T_{t-s_n}^*,
\qquad t>s_n.
 \eeqnn
Under this modification, $(\widetilde U_t)_{t>0}$ satisfies
$\widetilde U_{r+t} = T_t\widetilde U_r T_t^*$ for all $r,t>0$,
while (\ref{2.7}) remains unchanged. \qed

\btheorem\label{t2.2}
If $(\mu_t)_{t\ge0}$ is an SC-semigroup
with decomposition (\ref{2.1}), then we can write
 \beqlb\label{2.19}
\int_{H^\circ} K(x,a) M_t(dx)
=
\int_0^tds \int_{H^\circ} K(x,a) L_s(dx),
\quad t\ge0,a\in H,
 \eeqlb
where $L_s(dx)$ is a $\sigma$-finite kernel from $(0,\infty)$
to $H^\circ$ satisfying $L_{r+t} = (T_tL_r)|_{H^\circ}$ for all $r,t>0$ and
 \beqnn
\int_0^t ds \int_H (1\land\|x\|^2)L_s(dx) <\infty,
\qquad t\ge0.
 \eeqnn
\etheorem

\noindent{\it Proof.} If $(\mu_t)_{t\ge0}$ is an SC-semigroup given by
(\ref{2.1}), then $t\mapsto M_t$ is non-decreasing by the second equation in
(\ref{2.5}). Let $c\ge1$ and $b\ge0$ be as in the proof of Lemma~\ref{l2.1}
and let
 \beqnn
h(t) := \int_{H^\circ} (1\land\|x\|^2) M_t(dx),
\qquad t\ge0.
 \eeqnn
By (\ref{2.5}) we have, for $r,t\ge 0$,
 \beqnn
h(r+t) - h(r)
=
\int_{H^\circ} (1\land\|T_rx\|^2) M_t(dx),
 \eeqnn
which is bounded above by $c^2 e^{2br} h(t)$. As in the proof of
Lemma~\ref{l2.1}, one sees that $h(t)$ is absolutely continuous in
$t\ge0$. Since the family of finite measures $\nu_t(dx) :=
(1\land\|x\|^2) M_t(dx)$ is non-decreasing and $t\mapsto
h(t) = \nu_t(H^\circ)$ is absolutely continuous, $\nu([0,t],B) =
\nu_t(B)$ defines a locally bounded Borel measure $\nu(\cdot, B)$
on $[0,\infty)$ for each $B\in {\scr B}(H^\circ)$. A monotone
class argument shows that $\nu(A,\cdot)$ is a Borel measure on
$H^\circ$ for each $A\in {\scr B}([0,\infty))$, so that $\nu
(\cdot,\cdot)$ is a bimeasure. By \cite[p.502]{EK86}, there is a
probability kernel $J_s(dx)$ from $[0,\infty)$ to $H^\circ$ such
that
 \beqnn
\nu(A,B)
=
\int_A J_s(B)\nu(ds,H^\circ)
=
\int_A J_s(B)dh(s)
=
\int_A J_s(B)h^\prime(s) ds,
 \eeqnn
where $h^\prime(s)$ is the Radon-Nikodym derivative of $dh(s)$
relative to Lebesgue measure. Defining the $\sigma$-finite kernel
$L_s(dx) := (1\land \|x\|^2)^{-1} h^\prime(s) J_s(dx)$ we obtain
(\ref{2.19}). By the second equation of (\ref{2.5}) one can modify
the definition of $(L_t)_{t>0}$ so that $L_{r+t} =
(T_tL_r)|_{H^\circ}$ is satisfied for all $r,t>0$. \qed

We say the linear part $(b_t)_{t\ge0}$ of (\ref{2.1}) is {\it
absolutely continuous} if there exists an $H$-valued path
$(c_s)_{s>0}$ such that $\<b_t,a\> = \int_0^t\<c_s,a\>ds$ for all
$t\ge0$ and $a\in H$. The following theorem gives a Hilbert
space version of Li \cite[Theorem~2]{L95} or
\cite[Theorem~3.2]{L02} and extends the characterization of
Bogachev {\it et al} \cite[Lemma~2.6 and Proposition~4.3]{BRS96}.

\btheorem\label{t2.3} Suppose that $(\mu_t)_{t\ge0}$ is a family
of probability measures on $H$. If there is a family of infinitely
divisible probabilities $(\nu_s)_{s>0}$ such that $\nu_{r+t} =
T_t\nu_r$ for all $r,t>0$ and
 \beqlb\label{2.20}
\hat\mu_t(a)
=
\exp\bigg\{\int_0^t \log\hat\nu_s(a) ds\bigg\},
\qquad t\ge0, a\in H,
 \eeqlb
then $(\mu_t)_{t\ge0}$ is an SC-semigroup. Conversely, every
SC-semigroup $(\mu_t)_{t\ge0}$ with absolutely continuous linear
part has representation (\ref{2.20}). \etheorem

\noindent{\it Proof.} If $(\mu_t)_{t\ge0}$ is given by
(\ref{2.20}), it is clearly an SC-semigroup. Conversely, let
$(\mu_t)_{t\ge0}$ be an SC-semigroup and let $(U_s)_{s>0}$ and
$(L_s)_{s>0}$ be provided by Theorems~\ref{t2.1} and~\ref{t2.2}.
Suppose that $\<b_t,a\> = \int_0^t\<c_s,a\>ds$. By (\ref{2.6}), we
can modify the definition of $(c_s)_{s>0}$ so that
 \beqnn
c_{r+t}
=
T_tc_r
+ \int_{H^\circ} \left(\chi_{[0,1]}(\|T_tx\|)
- \chi_{[0,1]}(\|x\|)\right)T_tx L_r(dx),
\quad r,t>0.
 \eeqnn
Then we have the result by letting $\nu_s$ be the infinitely
divisible probability defined by the triple $(c_s,U_s,L_s)$. \qed

We may call the family $(\nu_s)_{s>0}$ in Theorem~\ref{t2.3} an
{\it entrance law} for $(T_t)_{t\ge0}$. (More precisely, it is an
entrance law for the deterministic Markov process $\{T_tx:
t\ge0\}$, as, for example, in Sharpe \cite{S88}.) If there is a probability
measure $\nu_0$ on $H$ such that $\nu_s = T_s\nu_0$ for all $s>0$,
we say that $(\nu_s)_{s>0}$ is {\it closable}. In this case, the
corresponding SC-semigroup $(\mu_t)_{t\ge0}$ is given by
 \beqlb\label{2.21}
\hat\mu_t(a)
=
\exp\bigg\{\int_0^t \log\hat\nu_0(T^*_sa) ds\bigg\},
\quad t\ge0, a\in H,
 \eeqlb
which belongs to the class (\ref{1.3}). This explains the
connection of our characterization with that of Bogachev {\it et
al} \cite{BRS96}.

Theorem~\ref{t2.3} gives a characterization for all SC-semigroups
under the assumption of absolute continuity on the linear part
$(b_t)_{t\ge0}$. This assumption cannot be removed since
$(b_t)_{t\ge0}$ can be discontinuous as pointed out in Schmuland
and Sun \cite{SS01}. The following example shows that it can even be
continuous but nowhere differentiable.

\medskip

{\bf Example 2.1} \ Consider $H=L^2([0,2\pi))$ and let $T_t$ be the
shift operator by $t\ge0$ (mod $2\pi$). For $t\ge0$ and $f\in
L^2([0,2\pi))$ set $b_t = (I-T_t)f$. Then $(\delta_{b_t})_{t\ge0}$
is a constant SC-semigroup. Taking the inner product against $f$
we obtain
 \beqnn
\<f,b_t\>
&=& \|f\|^2 - \frac{1}{2\pi}\int_0^{2\pi} f(x-t)f(x) dx  \\
&=& \|f\|^2 - |\hat f(0)|^2 - 2\sum_{n=1}^\infty |\hat f(n)|^2 \cos(nt),
 \eeqnn
where $\hat f$ is the Fourier transform of $f$. Now let $f$ be the
function whose Fourier coefficients are given by
 \beqnn
\hat f(n)
= \cases{2^{-k/2} &if $|n|=2^k$, $k\ge 1$,
\cr 0             &otherwise.}
 \eeqnn
Then we have
 \beqnn
\<f,b_t\> = 2- 2\sum_{k=1}^\infty 2^{-k}\cos(2^kt),
 \eeqnn
which is (up to a constant) Weierstrass's nowhere differentiable
continuous function.

\medskip

Let us consider another important special type of SC-semigroup
given by (\ref{2.1}) under the assumption:
 \beqlb\label{2.22}
\int_{H^\circ} (\|x\|\land\|x\|^2) M_t(dx)
<
\infty,
\quad t\ge 0.
 \eeqlb
Since (\ref{2.2}) holds automatically, (\ref{2.22}) is only a
first norm-moment condition on the restriction of $M_t$ to $\{x\in
H: \|x\| \ge1\}$. We say the SC-semigroup $(\mu_t)_{t\ge0}$ is
{\it centered} if
 \beqnn
\int_H \<x,a\> \mu_t(dx) = 0,
\quad t\ge 0, a\in H.
 \eeqnn
In this case, Theorem~\ref{t2.3} implies that
 \beqlb\label{2.23}
\hat\mu_t(a)
=
\exp\bigg\{\int_0^t\bigg[- \frac{1}{2}\<U_sa,a\>
+ \int_{H^\circ} K_1(x,a) L_s(dx)\bigg]ds\bigg\},
\qquad t\ge0, a\in H,
 \eeqlb
where
 \beqnn
K_1(x,a)
:=
e^{i\< x,a\>} -1 - i\< x,a\>,
\qquad x, a\in H.
 \eeqnn
Construction and regularity of OU-processes defined by
(\ref{2.23}) are discussed systematically in Dawson and Li
\cite{DL03}.

The characterizations (\ref{2.20}) and (\ref{2.23}) are of
interest since they include some SC-semigroups arising naturally
in applications which are not included in (\ref{1.3}) and
(\ref{2.21}). We shall see in the next two sections that a rich
class of such SC-semigroups arise in the study of fluctuation
limits of catalytic branching superprocesses with immigration. Two
particular examples are given below. We consider the Hilbert
space $L^2(0,\infty)$. Let
 \beqnn
g_t(x) = \frac{1}{\sqrt{2\pi t}} \exp\{-x^2/2t\},
\qquad t>0, x\in\IR
 \eeqnn
and
 \beqlb\label{2.24}
p_t(x,y) = g_t(x-y) - g_t(x+y),
\qquad t>0, \ x,y\in (0,\infty).
 \eeqlb
Then the transition semigroup $(P_t)_{t\ge0}$ of the ABM in
$(0,\infty)$ is defined by $P_0f =f$ and
 \beqlb\label{2.25}
P_tf(x) = \int_0^\infty p_t(x,y) f(y) dy,
\qquad t>0, x\in (0,\infty).
 \eeqlb
Let
 \beqlb\label{2.26}
k_t(y) = 2^{-1}(d/dx) p_t(x,y)|_{x=0^+} = yg_t(y)/t,
\qquad t>0, y\in (0,\infty).
 \eeqlb
It is not hard to check that
 \beqlb\label{2.27}
\int_0^\infty k_t(y)dt =1 \quad\mbox{and}\quad
k_{r+t}(y) = \int_0^\infty p_t(x,y)k_r(x)dx
 \eeqlb
for $r,t>0$ and $y\in (0,\infty)$.

\medskip

\noindent{\bf Example 2.2} \ Let $c>0$ and $x_0>0$. By
Theorem~\ref{t4.1}, there is a centered Gaussian SC-semigroup
$(\mu_t)_{t\ge0}$ on $L^2(0,\infty)$ given by
 \beqlb\label{2.28}
\hat \mu_t (f)
=
\exp\bigg\{-c\int_0^tP_sf(x_0)^2ds\bigg\},
\qquad t\ge0, f\in L^2(0,\infty).
 \eeqlb
This is a special form of (\ref{2.20}) and (\ref{2.23}) with
$(\nu_s)_{s>0}$ defined by
 \beqlb\label{2.29}
\hat \nu_s (f)
=
\exp\left\{-cP_sf(x_0)^2\right\},
\qquad s>0, f\in L^2(0,\infty).
 \eeqlb
Observe that $f\mapsto P_sf(x_0)^2$ is a well-defined functional
on $L^2(0,\infty)$ only for $s>0$. Thus, the SC-semigroup
(\ref{2.28}) is not included in (\ref{1.3}) and (\ref{2.21}).

\medskip
\noindent{\bf Example 2.3} \ Suppose that $(1\vee |u|)m(du)$ is a
finite measure on $\IR^\circ := \IR \setminus \{0\}$ and let
 \beqnn
\varphi(z)
=
\int_{\IR^\circ}\left(e^{iuz}-1-iuz\right)m(du),
\quad z\in \IR.
 \eeqnn
By Theorem~\ref{t4.3},
 \beqlb\label{2.30}
\hat \mu_t^\prime (f)
=
\exp\bigg\{\int_0^t\varphi(\<k_s,f\>)ds\bigg\},
\qquad t\ge0, f\in L^2(0,\infty)
 \eeqlb
defines a centered SC-semigroup $(\mu_t^\prime)_{t\ge0}$ on
$L^2(0,\infty)$. By (\ref{2.27}) one may check that
$(\mu_t^\prime)_{t\ge0}$ is included in (\ref{2.20}) and
(\ref{2.23}). Unless $m(\IR^\circ)=0$, this SC-semigroup is not
included in (\ref{1.3}) and (\ref{2.21}).

\medskip

\section{Fluctuation limits of superprocesses}
\setcounter{equation}{0}

In this section, we discuss small branching fluctuation limits of
catalytic branching superprocesses with immigration, which lead to
a class of OU-processes taking distribution values. Similar
fluctuation limits for superprocesses with function-valued
catalysts have been discussed in Gorostiza \cite{G96}, Gorostiza
and Li \cite{GL98, GL00}, and Li \cite{L99}. We shall only give an
outline of the arguments and refer the reader to the earlier
papers for details. As pointed out in \cite{L99}, the small
branching fluctuation limit is typically equivalent to the high
density and the large scale fluctuation limits. For simplicity, we
restrict to the case where the underlying motion is an ABM in $D
:= (0,\infty)$. We write $D$ instead of $(0,\infty)$ for the
underlying space in the sequel since $(0,\infty)$ and $[0,\infty)$
will appear frequently with quite different meanings. This
notation also suggests that some of the results can be modified to
the case where $D$ is a more general domain in $\IR^d$.

Let $M(D)$ denote the space of finite Borel measures on $D$
endowed with the topology of weak convergence. Let $\{B_t:
t\ge0\}$ be an ABM in $D$ with transition semigroup
$(P_t)_{t\ge0}$ defined by (\ref{2.25}). Let $\phi(\cdot,\cdot)$
be a function on $D\times [0,\infty)$ given by
 \beqlb\label{3.1}
\phi(x,z)
=
c(x)z^2 + \int_0^\infty(e^{-zu}-1+zu) m(x,d u),
\quad z\ge0,x\in D,
 \eeqlb
where $c\in B(D)^+$ and $u^2 m(x,du)$ is a bounded kernel from $D$
to $(0,\infty)$. For any $\eta\in B(D)^+$, there is a superprocess
in $M(D)$ with transition semigroup $(Q_t)_{t\ge0}$ determined by
 \beqlb\label{3.2}
\int_{M(D)} e^{-\nu(f)} Q_t(\mu,d\nu)
=
\exp\left\{-\mu(V_tf)\right\},
\qquad  f\in B(D)^+,
 \eeqlb
where $(V_t)_{t\ge0}$ is a semigroup of non-linear operators on
$B(D)^+$ defined by
 \beqlb\label{3.3}
V_tf(x) = P_tf(x)
- \int_0^td s\int_D \phi(y,V_sf(y)) \eta(y) P_{t-s}(x,dy),
\quad t\ge0, x\in D,
 \eeqlb
see for example Dawson \cite{D93}. The superprocess describes the
catalytic reaction of a large number of infinitesimal particles
moving according to the transition law of the ABM and splitting
according to the branching mechanism given by $\phi(\cdot,\cdot)$.
The value $\eta(x)$ represents the density at $x\in D$ of a
catalyst which causes the splitting. However, there are some
catalytic reactions in which the catalyst is concentrated on a
very small set and in that case the coefficient $\eta(\cdot)$ has
to be replaced by an irregular one, as in Pagliaro and Taylor
\cite{PT88}. These lead to the study of a catalyst given not by a
regular density function but rather by a measure $\eta\in M(D)$
with $\eta(dx):=$ ``catalytic mass in the volume element $dx$''.
Then we reformulate (\ref{3.3}) as
 \beqlb\label{3.4}
V_tf(x) = P_tf(x) -
\int_0^td s\int_D \phi(y,V_sf(y)) p_{t-s}(x,y)\eta(d y),
\qquad t\ge0, x\in D,
 \eeqlb
where $p_t(x,y)$ is given by (\ref{2.24}). A Markov process in
$M(D)$ with transition semigroup $(Q_t)_{t\ge0}$ defined by
(\ref{3.2}) and (\ref{3.4}) is called a {\it catalytic branching
super ABM} with parameters $(\eta,\phi)$.

Let $\{l_s(y): s>0,y>0\}$ be a continuous version of the local
time of $\{B_t: t\ge0\}$. Then $K(r,t) = \eta(l_t) - \eta(l_r)$
defines an additive functional of $\{B_t: t\ge0\}$. In view of
(\ref{2.24}) it is easy to check that
 \beqnn
\E_x\{K(0,t)\}
=
\int_0^tds\int_D p_s(x,y) \eta(dy)
\le
\eta(1)\int_0^t\frac{1}{\sqrt{2\pi s}}ds
\le
\eta(1)\sqrt{2t/\pi}.
 \eeqnn
Thus, $K(r,t)$ is admissible in the sense of \cite[p.49]{D94} and
the existence of the catalytic branching super ABM follows by
\cite[p.52]{D94}; see also \cite{Le00}. The study of
superprocesses with irregular catalysts was initiated by Dawson
and Fleischmann \cite{DF91, DF92} and there has been a
considerable development in the theory since then; see
Dawson and Fleischmann \cite{DF00} for a recent survey.

Set $\kappa_t(dx) = k_t(x)dx$. By (\ref{2.27}), $(\kappa_t)_{t>0}$
forms an entrance law for the underlying semigroup $(P_t)_{t\ge0}$,
that is, $\kappa_r P_t = \kappa_{r+t}$ for all $r,t>0$. Let
 \beqlb\label{3.5}
S_t(\kappa,f)
=
\kappa_t(f)
- \int_0^td s\int_D \phi(y,V_sf(y)) k_{t-s}(y)\eta(d y),
\quad t>0, f\in B(D)^+.
 \eeqlb
As in Li \cite{L96} one may see that
 \beqlb\label{3.6}
\int_{M(D)} e^{-\nu(f)} Q^\kappa_t(\mu,d\nu)
=
\exp\bigg\{-\mu(V_tf) - \int_0^t S_r(\kappa,f)d r\bigg\},
\quad f\in B(D)^+
 \eeqlb
defines the transition semigroup $(Q^\kappa_t)_{t\ge0}$ of a
Markov process $\{Y_t: t\ge0\}$ in $M(D)$, which we shall call a
{\it catalytic branching immigration super ABM with parameters}
$(\eta, \phi, \kappa)$. By (\ref{3.4}) and (\ref{3.5}) it is not
hard to check that
 \beqnn
\frac{d}{d\theta}V_t(\theta f)(x)\bigg|_{\theta=0^+}
=
P_tf(x) \quad\mbox{and}\quad
\frac{d}{d\theta}S_t(\kappa,\theta f)\bigg|_{\theta=0^+}
=
\kappa_t(f),
 \eeqnn
which, together with (\ref{3.6}), imply that
 \beqlb\label{3.6a}
\int_{M(D)} \nu(f) Q^\kappa_t(\mu,d\nu)
=
\mu(P_tf) + \int_0^t \kappa_r(f)d r.
 \eeqlb
By (\ref{2.27}) and (\ref{3.6a}) it follows that if $Y_0 =
\lambda$, then $\E\{Y_t(f)\} = \lambda(f)$ for all $t\ge0$ and
$f\in B(D)^+$, where $\lambda$ denotes Lebesgue measure.

Now we consider a small branching fluctuation limit of the
catalytic branching immigration ABM. For any $\theta>0$, let
$\phi_\theta(x,z) = \phi(x,\theta z)$ and ${\scr S}_\theta (D) =
\{\mu - \theta^{-1} \lambda: \mu\in M(D)\}$. Suppose that
$\{Y^\theta_t: t\ge0\}$ is a catalytic branching immigration ABM
with parameters $(\eta, \phi_\theta, \kappa)$ and $Y^\theta_0 =
\lambda$. As observed above, we have $\E\{Y^\theta_t(f)\} =
\lambda(f)$ for all $t\ge0$ and $f\in B(D)^+$. On the other hand,
$\phi_\theta(x,z) \to 0$ as $\theta \to 0$. By (\ref{2.4}),
(\ref{2.5}) and (\ref{2.6}) we have $Y^\theta_t(f) \to \lambda(f)$
in distribution as $\theta \to 0$. We define the fluctuation
process $\{Z_t^\theta: t\ge 0\}$ by
 \beqlb\label{3.7}
Z_t^\theta
=
\theta^{-1} [Y_t^\theta - \lambda],
\qquad t\ge 0.
 \eeqlb
As in Gorostiza and Li \cite{GL98} we see that $\{Z_t^\theta: t\ge
0\}$ is a centered signed-measure-valued Markov process with
transition semigroup $(R^\theta_t)_{t\ge0}$ determined by
 \beqnn
\int_{{\scr S}_\theta (D)} e^{-\nu(f)} R^\theta_t(\mu,d\nu)
=
\exp\bigg\{-\mu(\theta V^\theta_t(f/\theta))
+ \int_0^t\eta(\phi(\theta V^\theta_s(f/\theta)))ds\bigg\},
 \eeqnn
where $(V^\theta_t)_{t\ge0}$ is defined by
 \beqnn
V^\theta_tf(x)
+ \int_0^tds\int_D \phi_\theta(y,V^\theta_sf(y)) p_{t-s}(x,y)\eta(dy)
= P_tf(x).
 \eeqnn
Let ${\scr S}(D)$ be the space of infinitely differentiable
functions $f$ on $D$ such that
 \beqlb\label{3.8}
|||f|||_n
:=
\max_{0\le k\le n} \sup_{u\in D}
\bigg|(1+u^2)^n\frac{d^k}{du^k}f(u)\bigg|
< \infty,
\quad n=0,1,2,\dots.
 \eeqlb
Then ${\scr S}(D)$ topologized by the norms $\{|||\cdot|||_n:
n=0,1,2,\dots\}$ is a nuclear space. Let ${\scr S}^\prime(D)$
denote the dual space of ${\scr S}(D)$. As in \cite{GL98}, the
finite-dimensional distributions of $\{Z_t^\theta: t\ge 0\}$
converge as $\theta \to 0$ to those of an ${\scr
S}^\prime(D)$-valued Markov process $\{Z^0_t: t\ge 0\}$ with
transition semigroup $(R^0_t)_{t\ge0}$ determined by
 \beqlb\label{3.9}
\int_{{\scr S}^\prime(D)} e^{-\nu(f)} R^0_t(\mu,d\nu)
=
\exp\bigg\{-\mu(P_tf)
+ \int_0^t\eta(\phi(P_sf))ds\bigg\},
\quad f\in{\scr S}(D)^+.
 \eeqlb
Therefore, an OU-process with transition semigroup $(R^0_t)
_{t\ge0}$ is an approximation of the fluctuations of an
immigration process around the average.

\btheorem\label{t3.1} Let $\varphi(\cdot,\cdot)$ be a function on
$D \times\IR$ with the representation
 \beqlb\label{3.10}
\varphi(x,z) = -c(x)z^2 + \int_{\IR^\circ} (e^{izu}-1-izu) m(x,du),
\quad x\in D, z\in\IR,
 \eeqlb
where $c\in B(D)^+$ and $(|u|\land |u|^2) m(x,du)$ is a bounded
kernel from $D$ to $\IR^\circ := \IR \setminus \{0\}$. Then there
is a transition semigroup $(R_t)_{t\ge0}$ on ${\scr S}^\prime(D)$
given by
 \beqlb\label{3.11}
\int_{{\scr S}^\prime(D)} e^{i\nu(f)}R_t(\mu,d\nu)
=
\exp\bigg\{i\mu(P_tf) + \int_0^t\eta(\varphi(P_sf))ds\bigg\},
\quad t\ge0, f\in {\scr S}(D).
 \eeqlb
\etheorem

{\it Proof.} The semigroup $(R_t)_{t\ge0}$ can be obtained as in
\cite{GL98} by considering the difference of two Markov processes
with transition semigroups of the form (\ref{3.9}). \qed

Heuristically, an OU-process with transition semigroup $(R_t)
_{t\ge0}$ is the mixture of the fluctuations of two immigration
processes around their means. The branching mechanism of the
processes is determined by the function $\varphi (\cdot,\cdot)$
given by (\ref{3.10}) and the distribution of catalysts in $D$
that cause the branching is given by $\eta$. A more singular
transition semigroup is given by the following

\btheorem\label{t3.2}
Let $\varphi$ be a function on $\IR$ given by
 \beqlb\label{3.12}
\varphi(z)
=
- cz^2 + \int_{\IR^\circ}\left(e^{iuz}-1-iuz\right)m(du),
\qquad z\in \IR,
 \eeqlb
where $c\ge 0$ and $(|u|\land |u|^2) m(du)$ is a finite measure on
$\IR^\circ$. Then there is a transition semigroup $(R_t^\prime)
_{t\ge0}$ on ${\scr S}^\prime (D)$ given by
 \beqlb\label{3.13}
\int_{{\scr S}^\prime (D)} e^{i\nu(f)}R_t^\prime (\mu,d\nu)
=
\exp\bigg\{i\mu(P_tf) + \int_0^t \varphi(\kappa_s(f))ds\bigg\},
\quad t\ge0, f\in {\scr S}(D).
 \eeqlb
\etheorem

{\it Proof.} This transition semigroup is obtained from the one in
the last theorem by replacing $\varphi(x,z)$ and $\eta(dx)$ in
(\ref{3.11}) respectively by $\varphi(nz)$ and $\delta_{1/2n}
(dx)$ and letting $n\to \infty$. \qed

Roughly speaking, an OU-process with transition semigroup
$(R_t^\prime) _{t\ge0}$ represents the fluctuations of a process
over $D$ that branches very actively only near the absorbing
boundary.

\section{OU-processes with function values}
\setcounter{equation}{0}

In this section, we show that, under suitable conditions, the
OU-processes constructed in the last section take function values
from $L^2(D,\lambda)$.

Suppose that $\eta$ is a finite measure on $D$ and $\varphi
(\cdot,\cdot)$ is given by (\ref{3.10}) with $u^2m(x,du)$ being a
bounded kernel from $D$ to $\IR^\circ = \IR \setminus \{0\}$. Let
$W(ds,dx)$ be a white noise on $[0,\infty) \times D$ with
covariance measure $2c(x) ds\eta(dx)$ and $N(ds,du,dx)$ be a
Poisson random measure on $[0,\infty) \times \IR^\circ \times D$
with intensity $ds m(x,du) \eta(dx)$. Suppose that $W(ds,dx)$ and
$N(ds,du,dx)$ are defined on some complete probability space
$(\itOmega, {\scr F}, \P)$ and are independent of each other. Set
$\tilde N(ds,du,dx) = N(ds,du,dx) - ds m(x,du) \eta(dx)$. Then we
have

\btheorem\label{t4.1}
For each $t\ge0$, the function
 \beqlb\label{4.1}
Z^0_t(\omega,y)
:=
\int_0^t\int_D p_{t-s}(x,y) W(\omega,ds,dx)
+ \int_0^t\int_{\IR^\circ}\int_D up_{t-s}(x,y) \tilde N(\omega,ds,du,dx)
 \eeqlb
is well-defined in the $L^2(\itOmega \times D, \P \times \lambda)$
sense and $\{Z^0_t: t\ge0\}$ is a Markov process with state space
$L^2(D, \lambda)$, initial value zero and transition semigroup
$(R_t)_{t\ge0}$ given by
 \beqlb\label{4.2}
\int_{L^2(D, \lambda)} e^{i\<h,f\>}R_t(g,dh)
=
\exp\bigg\{i\<g,P_tf\> + \int_0^t\eta(\varphi(P_sf))ds\bigg\},
\quad f\in L^2(D, \lambda).
 \eeqlb
Moreover, $Z_t^0 (\omega,y)$ can be chosen as a function of $(t,
\omega, y)$ belonging to $L^2([0,T] \times \itOmega \times D,
\lambda \times \P \times \lambda)$ for each $T>0$. \etheorem

{\it Proof.} By the inequality
 \beqnn
\int_D p_{t-s}(x,y)^2 dy
<
\frac{1}{2\pi (t-s)}
\int_{\IR} \exp\bigg\{-\frac{y^2}{(t-s)}\bigg\} dy
=
\frac{1}{2\sqrt{\pi(t-s)}},
 \eeqnn
we have
 \beqnn
& &\int_D \E\bigg\{\bigg(\int_0^t\int_D p_{t-s}(x,y)
W(ds,dx) \bigg)^2\bigg\}dy  \\
&=&
2\int_D dy \int_0^tds \int_Dp_{t-s}(x,y)^2c(x)\eta(dx)  \\
&\le&
\int_0^t\frac{1}{\sqrt{\pi(t-s)}}ds \int_Dc(x)\eta(dx)  \\
&<& \infty
 \eeqnn
and
 \beqnn
& &\int_D \E\bigg\{\bigg(\int_0^t\int_{\IR^\circ}\int_D
up_{t-s}(x,y) \tilde N(ds,du,dx) \bigg)^2\bigg\}dy  \\
&=&
\int_D dy \int_0^tds \int_D\eta(dx) \int_{\IR^\circ}
u^2p_{t-s}(x,y)^2 m(x,du) \\
&\le&
\int_0^t\frac{1}{2\sqrt{\pi(t-s)}}ds \int_D\eta(dx)
\int_{\IR^\circ} u^2 m(x,du)  \\
&<& \infty.
 \eeqnn
Then the right hand side of (\ref{4.1}) is well-defined in the
$L^2(\itOmega \times D, \P \times \lambda)$ sense. By the same
reasoning, we see that it is also well-defined in the $L^2([0,
T] \times \itOmega \times D, \lambda \times \P \times
\lambda)$ sense. For any $f\in L^2(D,\lambda)$, we have
 \beqnn
\E\exp\bigg\{i\int_0^t\int_D P_{t-s}f(x) W(ds,dx)\bigg\}
=
\exp\bigg\{-\int_0^tds\int_D c(x)[P_{t-s}f(x)]^2 \eta(dx)\bigg\}
 \eeqnn
and
 \beqnn
& &\E\exp\bigg\{i\int_0^t\int_{\IR^\circ}\int_D uP_{t-s}f(x)
\tilde N(ds,du,dx)\bigg\} \\
&=&
\exp\bigg\{\int_0^tds \int_Dc(x)\eta(dx) \int_{\IR^\circ}
\left(\exp\{iuP_{t-s}f(x)\} - 1 - iuP_{t-s}f(x) \right)
m(x,du)\bigg\}.
 \eeqnn
Thus $\{Z^0_t: t\ge0\}$ has the asserted one-dimensional
distributions. If $g\in L^2(D,\lambda)$, then $P_tg\in
L^2(D,\lambda)$ for all $t\ge0$. Clearly, the distribution
$R_t(g,\cdot)$ of $P_tg+Z^0_t$ has characteristic functional given
by (\ref{4.2}) and $(R_t)_{t\ge0}$ is a transition semigroup on
$L^2(D,\lambda)$. The Markov property of $\{Z^0_t: t\ge0\}$
follows by a similar calculation of the characteristic functionals
of the finite-dimensional distributions. \qed

Suppose that $\varphi(\cdot)$ is given by (\ref{3.12}) with
$u^2m(du)$ being a finite measure on $\IR^\circ$. Set $\gamma(dx)
= (1-e^{-x^2})dx$ for $x\in D$. Let $\{B(t): t\ge0\}$ be a
one-dimensional Brownian motion with increasing process $2ct$ and
 $N(ds,du)$ be a Poisson random measure on $[0,\infty) \times
\IR^\circ$ with intensity $dsm(du)$. Suppose that $\{B(t):
t\ge0\}$ and $N(ds,du)$ are defined on some complete probability
space $(\itOmega, {\scr F}, \P)$ and are independent of each
other. Set $\tilde N(ds,du) = N(ds,du) - dsm(du)$. Then we have

\btheorem\label{t4.2}
For each $t\ge0$, the function
 \beqlb\label{4.3}
Z^0_t(\omega,y)
:=
\int_0^tk_{t-s}(y) B(\omega,ds)
+ \int_0^t\int_{\IR^\circ} uk_{t-s}(y) \tilde N(\omega,ds,du)
 \eeqlb
is well-defined in the $L^2(\itOmega \times D, \P \times \gamma)$
sense and $\{Z^0_t: t\ge0\}$ can be regarded as a Markov process
with state space ${\scr S}^\prime(\IR)$, initial value zero and
transition semigroup $(R_t^\prime)_{t\ge0}$ given by (\ref{3.13}).
Moreover, $Z_t^0 (\omega,y)$ can be chosen as a function of $(t,
\omega, y)$ belonging to $L^2([0,T] \times \itOmega \times D,
\lambda \times \P \times \gamma)$ for each $T>0$. \etheorem

{\it Proof.} For any $t>0$,
 \beqnn
\int_0^t k_s(y)^2 ds
\le
\int_0^\infty \frac{y^2}{2\pi s^3}e^{-y^2/s} ds
=
\frac{1}{2\pi y^2}.
 \eeqnn
Then we have
 \beqnn
& &\int_D \E\bigg\{\bigg(\int_0^t k_{t-s}(y)B(ds)\bigg)^2\bigg\}
 \gamma(dy)  \\
&=&
2c\int_D \gamma(dy) \int_0^tk_{t-s}(y)^2ds \\
&<& \infty
 \eeqnn
and
 \beqnn
& &\int_D \E\bigg\{\bigg(\int_0^t\int_{\IR^\circ} uk_{t-s}(y)
\tilde N(ds,du) \bigg)^2\bigg\}\gamma(dy)  \\
&=&
\int_D \gamma(dy) \int_0^t k_{t-s}(y)^2 ds \int_{\IR^\circ}
u^2 m(du) \\
&<& \infty.
 \eeqnn
Thus, the right hand side of (\ref{4.3}) is well-defined in the
$L^2(\itOmega \times D, \P \times \gamma)$ sense. Clearly, it is
also well-defined in the $L^2([0, T] \times \itOmega \times
D, \lambda \times \P \times \gamma)$ sense. For any $f\in L^2(D,
\lambda)$, we have
 \beqnn
\E\exp\bigg\{i\int_0^t \<k_{t-s},f\> B(ds)\bigg\}
=
\exp\bigg\{-\int_0^t c\<k_{t-s},f\>^2 ds\bigg\}
 \eeqnn
and
 \beqnn
& &\E\exp\bigg\{i\int_0^t\int_{\IR^\circ} u\<k_{t-s},f\>
\tilde N(ds,du)\bigg\} \\
&=&
\exp\bigg\{\int_0^tds \int_{\IR^\circ}
\left(\exp\{iu\<k_{t-s},f\>\} - 1
- iu\<k_{t-s},f\> \right) m(du)\bigg\}.
 \eeqnn
Therefore $\{Z^0_t: t\ge0\}$ has the correct one-dimensional
distributions. The asserted Markov property follows by a
calculation of the characteristic functionals of the
finite-dimensional distributions. \qed

\btheorem\label{t4.3}
If $(1\vee |u|) m(du)$ is a finite measure on $\IR^\circ$, then
for each $t\ge0$ the function
 \beqlb\label{4.4}
Z^0_t(y)
:=
\int_0^t\int_{\IR^\circ} uk_{t-s}(y) \tilde N(ds,du)
 \eeqlb
belongs to $L^2(D, \lambda)$ a.s.\ and $\{Z^0_t: t\ge0\}$ is a
Markov process with state space $L^2(D, \lambda)$, initial value
zero and transition semigroup $(R_t^\prime)_{t\ge0}$ given by
 \beqlb\label{4.5}
\int_{L^2(D,\lambda)} e^{i\<h,f\>}R_t^\prime (g,dh)
=
\exp\bigg\{i\<g,P_tf\>
+ \int_0^t \varphi(\<k_s,f\>)ds\bigg\},
\quad f\in L^2(D,\lambda),
 \eeqlb
where $\varphi$ is given by (\ref{3.12}) with $c=0$.
\etheorem

{\it Proof.} For any $t>0$, we have
 \beqnn
\int_0^t k_s(y) ds
=
\int_{y^2/2t}^\infty \frac{1}{\sqrt{\pi u}}e^{-u} du,
 \eeqnn
which is bounded in $y\ge 0$ and dominated by
 \beqnn
\int_{y^2/2t}^\infty \frac{1}{\sqrt{\pi}}e^{-u} du
=
\frac{1}{\sqrt{\pi}}e^{-y^2/2t}
 \eeqnn
for $y\ge\sqrt {2t}$. Therefore,
 \beqnn
\int_0^tds\int_{\IR^\circ} uk_{t-s} m(du)
 \eeqnn
belongs to $L^2(D,\lambda)$ under our assumption. Since $k_t\in
L^2(D,\lambda)$ for every $t>0$ and a.s.\
 \beqnn
\int_0^t\int_{\IR^\circ} uk_{t-s} N(ds,du)
 \eeqnn
is a finite sum, we have $Z^0_t\in L^2(D, \lambda)$ a.s.\ If $g\in
L^2(D,\lambda)$, then $P_tg\in L^2(D,\lambda)$ for all $t\ge0$ and
the distribution $R_t(g,\cdot)$ of $P_tg+Z^0_t$ has characteristic
functional given by (\ref{4.5}). Clearly, $(R_t)_{t\ge0}$ is a
transition semigroup on $L^2(D,\lambda)$. The Markov property of
$\{Z^0_t: t\ge0\}$ follows by a calculation of the characteristic
functionals of the finite-dimensional distributions. \qed

As in Li \cite{L96} one may see that the generalized OU-processes
given by (\ref{4.2}) and (\ref{4.5}) usually do not have right
continuous sample paths, neither do they have the strong Markov
property. We shall prove in the next section that they do have
those properties if we regard them as processes in another
suitably chosen state space.

\section{OU-processes with signed-measure values}
\setcounter{equation}{0}

In this section, we show that some of the generalized OU-processes
given by (\ref{4.2}) and (\ref{4.5}) behave very regularly in the
space of signed-measures. Indeed, from the proof of
Theorem~\ref{t5.1} we know that they are essentially special forms
of the immigration processes studied in Li \cite{L95, L96}.

Given a locally compact metric space $E$, we denote by $M(E)$ the
space of finite Borel measures on $E$. Let $\{f_n\}_{n=1}^\infty$
be a dense subset of the space of all bounded uniformly continuous
functions on $E$. We define the metric $r(\cdot,\cdot)$ on
$M(E)$ by
 \beqlb\label{5.1}
r(\mu,\nu)
=
\sum_{n=1}^\infty 2^{-n}(1\land |\mu(f_n) - \nu(f_n)|),
\qquad \mu,\nu\in M(E).
 \eeqlb
Clearly, this metric is compatible with the topology of weak
convergence in $M(E)$. Let $S(E) = \{\mu^+-\mu^-: \mu^+, \mu^- \in
M(E)\}$ be the space of finite signed-measures on $E$. Define a
metric $\rho(\cdot,\cdot)$ on $S(E)$ by
 \beqlb\label{5.2}
\rho(\mu,\nu)
&=&
\inf\{r(\mu^+,\nu^+) + r(\mu^-,\nu^-):
\mu^+,\mu^-,\nu^+,\nu^-\in M(E)  \nonumber \\
& &
\qquad\qquad \mbox{ with }\mu^+ - \mu^-
= \mu
\mbox{ and } \nu^+-\nu^- = \nu\}.
 \eeqlb
Then $\mu_n \to \mu_0$ in $S(E)$ if and only if there are
decompositions $\mu_n = \mu^+_n - \mu^-_n$ and $\mu_0 = \mu^+_0 -
\mu^-_0$ such that $\mu^+_n \to \mu^+_0$ and $\mu^-_n\to \mu^-_0$
in $M(E)$. Below, we shall consider the metric space $(S(E),\rho)$
for $E= (0,\infty)$ or $[0,\infty)$.

Suppose that $\eta$ is a finite measure, $\varphi(\cdot,\cdot)$
is given by (\ref{3.10}) with $c(x) \equiv 0$, $|u|
m(x,du)$ is a bounded kernel from $D$ to $\IR^\circ$, and
that $N(ds,du,dx)$ is a Poisson random measure on $[0,\infty)
\times \IR^\circ \times D$ with intensity $ds m(x,du) \eta(dx)$.
Let $\tilde N(ds,du,dx) = N(ds,du,dx) - ds m(x,du) \eta(dx)$ and
 \beqlb\label{5.3}
Y_t(f)
:=
\int_0^t\int_{\IR^\circ}\int_D uP_{t-s}f(x) \tilde N(ds,du,dx),
\qquad t\ge0, f\in B(D).
 \eeqlb

\btheorem\label{t5.1} The process $\{Y_t: t\ge0\}$ defined by
(\ref{5.3}) is an a.s.\ right continuous $S(D)$-valued strong
Markov process with transition semigroup $(R_t)_{t\ge0}$ defined
by
 \beqlb\label{5.4}
\int_{S(D)} e^{i\nu(f)}R_t(\mu,d\nu)
=
\exp\bigg\{i\mu(P_tf) + \int_0^t\eta(\varphi(P_sf))ds\bigg\},
\quad f\in B(D).
 \eeqlb
\etheorem

{\it Proof.} We define the positive part $\{Y^+_t: t\ge0\}$ of
$\{Y_t: t\ge0\}$ by

 \beqnn
Y^+_t(f)
:=
\int_0^t\int_0^\infty\int_D uP_{t-s}f(x) N(ds,du,dx),
\qquad t\ge0, f\in B(D).
 \eeqnn
By the assumptions,
 \beqnn
\E\{Y^+_t(f)\}
&=&
\int_0^tds\int_D \eta(dx) \int_0^\infty uP_{t-s}f(x) m(x,du)   \\
&\le&
t\|f\|\eta(D) \sup_{x\in D} \int_0^\infty u m(x,du)   \\
&<&
\infty.
 \eeqnn
Then $\{Y^+_t: t\ge0\}$ is a well-defined $M(D)$-valued process,
which is clearly a special form of the immigration process
considered in \cite{L96} without branching. By \cite[Theorem~4.1]{L96},
$\{Y^+_t: t\ge0\}$ is a.s.\ right continuous. Similarly, the
negative part $\{Y^-_t: t\ge0\}$ of $\{Y_t: t\ge0\}$ defined by
 \beqnn
Y^-_t(f)
:=
- \int_0^t\int_{-\infty}^0\int_D uP_{t-s}f(x) N(ds,du,dx),
\qquad t\ge0, f\in B(D)
 \eeqnn
is also an a.s.\ right continuous immigration process. Then one can
easily see that $\{Y_t: t\ge0\}$ defined by (\ref{5.3}) is an a.s.\
right continuous $S(D)$-valued Markov process with transition
semigroup $(R_t)_{t\ge0}$. The strong Markov property holds since
$(R_t)_{t\ge0}$ is clearly Feller. \qed

Suppose that $\varphi$ is given by (\ref{3.12}) with $c=0$,
 $|u| m(du)$ is a finite measure on $\IR^\circ$, and
  $N(ds,du)$ is a Poisson random measure on $[0,\infty)
\times \IR^\circ$ with intensity $dsm(du)$. Let $\tilde N(ds,du) =
N(ds,du) - dsm(du)$, and
 \beqlb\label{5.5}
Y_t(f)
:=
\int_0^t\int_{\IR^\circ} u\kappa_{t-s}(f) \tilde N(ds,du),
\qquad t\ge0, f\in B(D).
 \eeqlb
By an argument similar to that in the proof of Theorem~\ref{t5.1}
we get

\btheorem\label{t5.2} The process $\{Y_t: t\ge0\}$ defined by
(\ref{5.5}) is an $S(D)$-valued Markov process with transition
semigroup $(R_t^\prime)_{t\ge0}$ defined by
 \beqlb\label{5.6}
\int_{S(D)} e^{i\nu(f)}R_t^\prime (\mu,d\nu)
=
\exp\bigg\{i\mu(P_tf) + \int_0^t \varphi(\kappa_s(f))ds\bigg\},
\quad f\in B(D).
 \eeqlb
\etheorem

As in \cite{L96} one can see that the process (\ref{5.5}) does not
have any right continuous modification. Observe that $h(x) :=
(1-e^{-x})$ is an excessive function of $(P_t)_{t\ge0}$ and
 \beqlb\label{5.7}
T_tf(x)
=
\left\{\begin{array}{ll}
h(x)^{-1}P_t(hf)(x)
&\mbox{ for $t>0$ and $x>0$,}     \\
2\kappa_t(hf) = (d/dx)P_t(hf)(0^+)
&\mbox{ for $t>0$ and $x=0$}
\end{array}\right.
 \eeqlb
defines the transition semigroup $(T_t)_{t\ge0}$ of a Markov
process on $[0,\infty)$.

\btheorem\label{t5.3} Let $\{Y_t: t\ge0\}$ be defined by
(\ref{5.5}), $Z_t(\{0\}) = 0$, and $Z_t(dx) =
(1-e^{-x})Y_t(dx)$ for $x>0$. Then $\{Z_t: t\ge0\}$ is an
$S([0,\infty))$-valued Markov process with transition semigroup
$(S_t)_{t\ge0}$ defined by
 \beqlb\label{5.8}
\int_{S([0,\infty))} e^{i\nu(f)} S_t(\mu,d\nu)
=
\exp\bigg\{i\mu(T_tf) + \int_0^t \varphi(\kappa_s(hf))ds\bigg\},
\quad f\in B([0,\infty)).
 \eeqlb
Moreover, $\{Z_t: t\ge0\}$ has a right continuous strong Markov
realization.
\etheorem

{\it Proof.} The first assertion holds by Theorem~\ref{t5.2}.
Observe that
 \beqnn
\int_{S([0,\infty))} e^{i\nu(f)} S_t(\mu,d\nu)
=
\exp\bigg\{i\mu(T_tf) + \int_0^t \varphi(2^{-1}T_sf(0))ds\bigg\},
\quad f\in B([0,\infty)),
 \eeqnn
by (\ref{5.7}) and (\ref{5.8}). Then the second assertion follows
from \cite[Theorem~4.1]{L96} as in the proof of
Theorem~\ref{t5.1}. \qed

\end{document}